\input amstex\documentstyle{amsppt}  
\pagewidth{12.5cm}\pageheight{19cm}\magnification\magstep1
\topmatter
\title Sharpness versus bluntness in affine Weyl groups\endtitle
\author G. Lusztig\endauthor
\address{Department of Mathematics, M.I.T., Cambridge, MA 02139}\endaddress
\thanks{Supported by NSF grant DMS-2153741}\endthanks
\endtopmatter   
\document

\define\Irr{\text{\rm Irr}}
\define\sneq{\subsetneqq}

\define\mpb{\medpagebreak}

\define\hS{\hat S}
\define\hW{\hat W}

\define\si{\sim}

\define\sqc{\sqcup}

\define\lb{\linebreak}

\define\part{\partial}
\define\emp{\emptyset}

\define\iy{\infty}
\define\m{\mapsto}
\define\do{\dots}

\define\sub{\subset}    

\define\T{\times}
\define\ti{\tilde}
\define\nl{\newline}
\redefine\i{^{-1}}

\define\un{\underline}

\define\bbq{\bar{\QQ}_l}

\define\a{\alpha}
\redefine\b{\beta}
\redefine\c{\chi}

\define\e{\epsilon}

\define\io{\iota}
\redefine\o{\omega}

\redefine\t{\tau}

\define\z{\zeta}
\define\x{\xi}

\redefine\D{\Delta}
\define\Om{\Omega}
\define\Si{\Sigma}
\define\Th{\Theta}

\define\CC{\bold C}

\define\NN{\bold N}

\define\QQ{\bold Q}

\define\SS{\bold S}

\define\ZZ{\bold Z}

\define\ca{\Cal A}
\define\cb{\Cal B}
\define\cc{\Cal C}
\define\cd{\Cal D}
\define\ce{\Cal E}
\define\cf{\Cal F}
\define\cg{\Cal G}

\define\fA{\frak A}

\define\fS{\frak S}

\define\tS{\ti S}

\define\tU{\ti U}

\define\tW{\ti W}

\head  Introduction\endhead
\subhead 0.1\endsubhead
Let $W$ be a Weyl group with set of simple reflections $S$ and let
$\o$ be an automorphism of $W$ preserving $S$ (in the non simply
laced case we assume in addition that $\o$ is as in 1.3).
Let $G(F_q)$ be an adjoint reductive group over the finite field $F_q$
with Weyl group $W$ and with twisting described by $\o$. Let
$Cu(W,\o)$ be the set of (isomorphism classes of) irreducible
unipotent cuspidal representations (over $\bbq$) of $G(F_q)$.
(We can regard this set as being independent of $F_q$.)
We are interested in three questions:

(i) Under what conditions is $Cu(W,\o)$ nonempty?

(ii) More precisely, how many elements does $Cu(W,\o)$ have?

(iii) How to classify the cuspidal local systems on unipotent classes
in a simply connected reductive group over $\CC$ with Weyl group
$W$? 

The answer to (i),(ii) is known from \cite{L84a}; it involves
among other things considerations of intersection cohomology of
certain algebraic varieties. The answer to (iii) is known from
\cite{L84b}; it involves considerations of perverse sheaves.
(It actually depends not only on $W$ but also on a root datum
which gives rise to $W$.)

We would like to find a way to predict the answers to (i),(ii),(iii)
which is not based on the considerations above but instead it is
based on pure algebra or more precisely on

(a) the knowledge of the generic degrees of
representations of the Iwahori-Hecke algebra of $W$.
\nl
For (i) this is done in 2.2 where the notion of sharpness of
$(W,\o)$ is defined by requiring that a certain quantity defined in
terms of (a) should reach its maximum. (See also 4.1.)

In this paper we will treat (ii) and (iii) as parts of another
problem.

Let $\cg$ be a split simple adjoint group over a $p$-adic field.
In \cite{L83} we defined the notion of unipotent representation
of $\cg$ and stated a refinement of a conjecture of Deligne and
Langlands on the classification of such representations. The conjecture
(and its extension to the not necessarily split case) was proved in
\cite{L95},\cite{L02} (a special case was proved in \cite{KL87}).
The method of \cite{L95},\cite{L02} involve among other things
defining two
collections of diagrams associated to an irreducible affine Weyl group
$W$ with an automorphism $\o$; one of them (said to be the arithmetic
diagrams)
is defined in terms of $\cg$ and is used to parametrize the unipotent
representations of $\cg$; the other one (said to be the geometric
diagrams) is defined in terms of the geometry of the dual group over
$\CC$. In {\it loc.cit.} it is shown that the arithmetic diagrams are in
bijection with the geometric diagrams, but the bijection was given in
each case without explanation. (In \cite{L95},\cite{L02}, the
arithmetic/geometric
diagrams appear with certain repetitions which
are ignored here.)

In this paper we restrict ourselves to the case of ``cuspidal''
arithmetic/geometric diagrams that is, arithmetic/geometric diagrams
corresponding to unipotent
supercuspidal representations of a simple $p$-adic group. 
We will show that the cuspidal arithmetic diagrams, the cuspidal
geometric
diagrams and the bijection between them can be defined by pure algebra (see (a)). (Very similar results hold also for not necessarily
cuspidal arithmetic/geometric diagrams, as will be shown elsewhere.) 

\subhead 0.2\endsubhead
Our approach is as follows. Using the notion of sharpness of $(W,\o)$
we define (see 2.4) the notion of sharp triples
associated to the affine Weyl group $W$. Then we restrict the notion
of sharpness to one of strict sharpness (see 2.6)
by requiring that a certain quantity reaches its
minimum. From this we obtain the notion of blunt pair.
(At this point one obtains a solution of 0.1(iii) by pure algebra,
see 4.2.)
From this we can recover the cuspidal geometric diagrams.
Now the sharp triples give rise to the cuspidal arithmetic 
diagrams (for classical types) while for exceptional types
there is an obvious way to enhance the sharp triples to get the
cuspidal arithmetic 
diagrams.
At this point one obtains a solution of 0.1(ii) by pure algebra,
see 4.3.

\subhead 0.3\endsubhead
In this paper $|?|$ denotes the cardinal of a finite set $?$
or the absolute value of a real number $?$.

\subhead 0.4\endsubhead
I wish to thank Zhiwei Yun for useful comments.

\head Contents\endhead
1. Permissible and admissible automorphisms of an affine Weyl group.

2. Sharpness.

3. Bluntness.

4. Comments.

\head 1. Permissible and admissible automorphisms of an affine Weyl
group\endhead
\subhead 1.1\endsubhead
Let $W$ be a Coxeter group with (finite) set $S$ of simple reflections;
let $l:W@>>>\NN$ be the length function.
Let $M(W)$ be the set of $2$-element subsets $\{s,s'\}\sub S$ such
that the order of $ss'$ (or $s's$) is finite and $>3$.
We say that $W$ is simply laced if $M(W)=\emp$.
An {\it orientation} of $W$ is a rule which to any
$\{s,s'\}\in M(W)$ associates a total order $<$ on $\{s,s'\}$.
If an orientation $E$ of $W$ is given we can define a new
orientation $E^*$ of $W$ (said to be opposed to the original one)
by replacing the total
order on each $\{s,s'\}\in M(W)$ by the opposite one.

We say that $W$ is oriented if an orientation of $W$ is specified.

Let $\Om=\Om_W$ be the group of automorphisms of $W$ which map $S$ onto
itself. 
For $\o\in\Om$ we denote by $S_\o$ the set of orbits of $\o:S@>>>S$.

Let $\o\in\Om$ be such that for any $x\in S_\o$, the subgroup $W_x$ of
$W$ generated by $x$ is finite;
we denote by $w_x$ the longest element of this subgroup.
For $x\ne x'$ in $S_\o$ let $n_{x,x'}\in\ZZ_{>0}\cup\{\iy\}$ be the
order of $w_xw_{x'}\in W$; we have $n_{x,x'}=n_{x',x}$. Let $\un W$ be
the Coxeter group with generators $\{x;x\in S_\o\}$ and relations
$x^2=1$ for $x\in S_\o$ and $(xx')^{n_{x,x'}}=1$ for any $x\ne x'$ in
$S_\o$ such that $n_{x,x'}<\iy$. Clearly, we have a unique group
homomorphism $\c:\un W@>>>W^\o:=\{w\in W;\o(w)=w\}$ such that
$\c(x)=w_x$ for all $x\in S_x$. From \cite{L14, Appendix, Theorem A8}
we see that

(a) {\it $\c:\un W@>>>W^\o$ is an isomorphism.}

\mpb

Thus, $W^\o$ is a Coxeter group with simple reflections
$\{w_x;x\in S_\o\}$.

We say that $\o\in\Om$ is {\it good} if for any $\{x,x'\}\in M(W^\o)$
we have $l(w_x)\ne l(w_{x'})$; in this case we can order $\{x,x'\}$ by
writing $x>x'$ if $l(w_x)<l(w_{x'})$; thus $W^\o$ is oriented when $\o$
is good.

\subhead 1.2\endsubhead
For $\o\in\Om$ we set $S^\o=\{s\in S;\o(s)=s\}$. Assuming that $W$ is
simply laced, we say that $\o\in\Om$ is {\it admissible} if 

(a) $S^\o\ne\emp$ and

(b) for any $s,s'$ in the same $\o$-orbit on $S$ we have $ss'=s's$.

\subhead 1.3\endsubhead
We now assume that $W,S$ is a Weyl group. Let $\fA_W$ be the set
of all $\o\in\Om_W$ such that the intersection of any
$\{s,s'\}\in M(W)$ with any $\o$-orbit in $S$ contains at most one
element. Note that if $\o\in\fA_W$ then any power of $\o$ is in
$\fA_W$.
For $\o\in\fA_W$ let $r(\o)$ be the number of $\o$-orbits on $S$;
let $ord(\o)$ be the order of $\o$. Let $op\in\fA_W$ be given
by conjugation by the longest element of $W$.

Given $\o\in\fA_W$ we say that $W$ is $\o$-irreducible
if $W$ is a product of irreducible Weyl groups which are
cyclically permuted by $\o$. Given $\o\in\fA_W$, one can write
canonically $W$ as a product of $\o$-stable, $\o$-irreducible
Weyl groups, said to be the $\o$-irreducible factors of $W$.

We now assume that $W$ is an irreducible Weyl group.
It has type $A_n (n\ge1)$, $B_n=C_n (n\ge2)$,
$D_n (n\ge4)$, $E_n (n=6,7,8)$, $F_4$ or $G_2$. (We then write
$W=A_n$ or $W=B_n$, etc.)
In this case $\fA_W$ is a subgroup of $\Om_W$.

We have $\fA_W=\{1\}$ if $W$ is $A_1$ or $B_n,n\ge2$ or
$E_7,E_8,F_4,G_2$.
We have $|\fA_W|=2$ if $W$ is $A_n,n\ge2$ or $D_n,n\ge5$ or $E_6$; let
$\o_2$ be the element of order $2$ in $\fA_W$.
If $W=D_4$, we have $\fA_W\cong\fS_3$ (symmetric group); let
$\o_2$ (resp. $\o_3$) be an element of order $2$ (resp. $3$) in
$\fA_W$.

Let $Fin_{or}$ be the set of (isomorphism classes of) 
irreducible oriented Weyl groups.
Let $\SS^{adm}_{fin}$ be the set of (isomorphism classes of)
pairs $(W,\o)$ where $W$ is an irreducible Weyl group of simply
laced type and $\o\in\fA_W=\Om_W$ is admissible.
We give below a list of the objects $\SS^{adm}_{fin}$; in each case
we describe also $W^\o$ (an object of $Fin_{or}$, since $\o$ is good).

$(A_n,1), n\ge1$, $W^\o=A_n$

$(D_n,1), n\ge4$, $W^\o=D_n$

$(E_n,1), n=6,7,8$, $W^\o=E_n$

$(A_{2n-1},\o_2),n\ge2$, $W^\o=C_n$

$(D_{n+1},\o_2),n\ge3$, $W^\o=B_n$

$(E_6,\o_2)$, $W^\o=F_4$

$(D_4,\o_3)$, $W^\o=G_2$.

Note that $(W,\o)\m W^\o$ defines a bijection

(a) $\SS^{adm}_{fin}@>>>Fin_{or}$.
\nl
Note also that for $n\ge3$, the Weyl group $B_n=C_n$ has two distinct
orientations (opposed to each other) which are denoted by $B_n,C_n$.

Now $F\m F^*$ (see 1.1) is an involution of $Fin_{or}$ which
interchanges $B_n,C_n$ for $n\ge3$ and is the identity on all
other objects.

\subhead 1.4\endsubhead
In the remainder of this section we assume that $W$ is an irreducible
affine Weyl group. It has type

$\ca_n (n\ge1)$, $\cb_n (n\ge3)$,
$\cc_n (n\ge2)$, $\cd_n (n\ge4)$, $\ce_n (n=6,7,8)$, $\cf_4$ or $\cg_2$. \nl
(We then write 
$W=\ca_n$ or $W=\cb_n$, etc.) We sometimes write $\cd_3$ instead
of $\ca_3$ and viceversa.

For any $I\sneq S$ we denote by $W_I$ the subgroup of $W$
generated by $I$; this is a (finite) Weyl group.

Let $S_*$ be the set of all $s\in S$ where the function
$s\m|W_{S-\{s\}}|$ from $S$ to $\ZZ_{>0}$ reaches its maximum value.

Let $T$ be the union of all finite conjugacy classes of $W$; this is a
free abelian group which is normal of finite index in $W$. 
Let $\Om^1=\Om_W^1$ be the group consisting of all $\o\in\Om$  
such that the restriction of $\o$ to $T$
is conjugation by some element of $W$. This is an abelian normal
subgroup of $\Om$. Now the conjugation action of
$\Om$ on $S$ preserves $S_*$; this restricts to a free
transitive action of $\Om^1$ on $S_*$.

\subhead 1.5\endsubhead
For any $d\ge1$ let $\Om^d=\Om^d_W$ be the set of all $\o\in\Om=\Om_W$
such that that the image of $\o$ in $\Om/\Om^1$ has order $d$. This
agrees with our earlier notation for $\Om^1$. We have
$\Om=\sqc_{d\ge1}\Om^d$. Note that $\Om^d$ is a union of conjugacy
classes in $\Om$.

If $W$ is $\ca_1$, then $\Om=\Om^1$ is cyclic of order 
$2$; we denote by $\o_2$ the generator of $\Om^1$.

If $W$ is $\ca_n,n\ge2$ then $\Om=\Om^1\sqc\Om^2$ is a dihedral
group of order $2(n+1)$ and $\Om^1$ is cyclic of order $n+1$;   
we denote by $\o_{n+1}$ a generator of $\Om^1$. If $n$ is even we
denote by $\o_2$ any element of $\Om^2$. If $n$ is odd we denote by 
$\o_2$ (resp. $\o'_2$) any element of $\Om^2$ such that
$S^{\o_2}\ne\emp$ (resp. $S^{\o'_2}=\emp$).

If $W$ is $\cb_n,n\ge3$ or $\cc_n,n\ge2$, then $\Om=\Om^1$ is  
cyclic of order $2$; we denote by $\o_2$ the generator of $\Om$.

If $W$ is $\cd_n,n\ge5$, then $\Om=\Om^1\sqc\Om^2$ is dihedral of order
$8$; $\Om^1$ is cyclic of order $4$ if $n$ is odd and non-cyclic of
order $4$ if $n$ is even. Let ${}'\Om=\{\o\in\Om;|S^\o|\le1\}$.
For $d=1,2$ we set ${}'\Om^d={}'\Om\cap\Om^d$.
We denote by $\o_{2,2}$ the element of order $2$ of $\Om^1$ such that
$|S^{\o_{2,2}}|=n-3$. 
We denote by $\o'_{2,2}$ an element of order $2$ in ${}'\Om$.
(There are two such elements. If $n$ is even they form ${}'\Om^1$. 
If $n$ is odd they form ${}'\Om^2$.) We denote by $\o_4$ an element of 
order $4$ of $\Om$; we have $\o_4\in{}'\Om$ (There are two such elements. 
If $n$ is odd, they form ${}'\Om^1$. If $n$ is even they form ${}'\Om^2$.)
We denote by $\o_2$ an element of order $2$ of $\Om^2$ such that
$|S^{\o_2}|=n-1$. (There are two such elements.) 

If $W$ is $\cd_4$, then $\Om$ is the symmetric group in
four letters and $\Om^1$ is non-cyclic of order $4$.
We denote by $\o_{2,2}=\o'_{2,2}$ an element of order $2$ of
$\Om^1$ such that $|S^{\o_{2,2}}|=1$. (There are three such
elements. They form together with $1$ the group $\Om^1$.) 
We denote by $\o_2$ an element of order $2$ of $\Om$
such that $|S^{\o_2}|=3$. (There are six such elements; they
are in $\Om^2$.) We denote by $\o_4$ an element of order $4$ of $\Om$. 
(There are $6$ such elements; they are in $\Om^2$.)
We denote by $\o_3$ an element of
order $3$ of $\Om$. (There are $8$ such elements. They form $\Om^3$.)

If $W$ is $\ce_6$, then $\Om=\Om^1\sqc\Om^2$ is the symmetric group
in three letters and $\Om^1$ is cyclic of order $3$; we denote by
$\o_3$ a generator of $\Om^1$ and by $\o_2$ an element of order
$2$ of $\Om$. (We have $\o_2\in\Om^2$.)

If $W$ is $\ce_7$, then $\Om=\Om^1$ is cyclic of order $2$;
we denote by $\o_2$ the generator of $\Om^1$.

If $W$ is $\ce_8,\cf_4$ or $\cg_2$, then $\Om=\Om^1=\{1\}$.

\subhead 1.6\endsubhead
For $s\in S$ we denote by $j(s)$ the smallest integer $j\ge0$
such that there exists a sequence $s_0,s_1,\do,s_j$ in $S$
with $s_0s_1\ne s_1s_0,s_1s_2\ne s_2s_1$, $\do$,
$s_{j-1}s_j\ne s_js_{j-1},s_0\in S_*,s_j=s$. We show:

(a) {\it $W$ has a canonical orientation.}
\nl
If $W$ is $\ca_n (n\ge1)$, $\cd_n (n\ge4)$, or $\ce_n (n=6,7,8)$,
this is obvious since $M(W)=\emp$.

If $W$ is  $\cb_n (n\ge3)$ or $\cf_4$, $\cg_2$ then $M(W)$ consists of
a single pair $\{s,s'\}$. We have $j(s)\ne j(s')$ and we say that
$s>s'$ if $j(s)<j(s')$ (this is the canonical orientation).

If $W$ is $\cc_n (n\ge2)$ then $M(W)$ consists of two pairs
$\{s_1,s'_1\}$, $\{s_2,s'_2\}$. We have $\{j(s_1),j(s'_1)\}=\{0,1\}$,
$\{j(s_2),j(s'_2)\}=\{0,1\}$. We say that $s_1>s'_1$ if $j(s_1)=0$ (that is,
$s_1\in S_*$); we say that $s_2>s'_2$ if $j(s_2)=0$ (that is,
$s_2\in S_*$); this is the canonical orientation.

Let $Af_{or}$ be the set of (isomorphism classes of) irreducible
oriented affine Weyl groups. Let $Af_{or}^{can}$ be the subset of
$Af_{or}$ consisting of irreducible affine Weyl groups with the
canonical orientation.

We have a bijection $h:Fin_{or}@>\si>>Af_{or}^{can}$; its inverse
takes an oriented $(W,S)$ in $Af_{or}^{can}$ to
$W_{S-\{s\}}$ (where $s\in S_*$) with the induced orientation.

If $W$ is $\cf_4$ or $\cg_2$, there are two orientations of $W$; the
canonical orientation (denoted again by $\cf_4$ or $\cg_2$) and
its opposite (denoted by $\cf_4^*$ or $\cg_2^*$).

If $W$ is $\cb_n (n\ge3)$, there are two orientations of $W$; the
canonical orientation is denoted by $E^{DB}_n$; its opposite is denoted
by $E^{DC}_n$.

If $W$ is $\cc_n (n\ge2)$, there are three orientations of $W$; the
canonical orientation is denoted by $E^{CC}_n$; its opposite is denoted
by $E^{BB}_n$. The third orientation is denoted by $E^{BC}_n$; it is
characterized by the property that it is isomorphic to its opposed one.

We also define $E^{DC}_2$ to be the same as $E^{BB}_2$.
The following is a complete list of objects of $Af_{or}$.

$\ca_n (n\ge1)$, $\cd_n,(n\ge4)$, $\ce_n (n=6,7,8)$, $\cf_4,\cf_4^*$,
$\cg_2,\cg_2^*$, $E_n^{DB},E_n^{DC}, (n\ge3)$,
$E_n^{CC},E_n^{BB}, E_n^{BC}, (n\ge2)$.
\nl
There is a unique involution $\t:Af_{or}@>>>Af_{or}$
such that $\t(h(F))=h(F^*)$, $\t(h(F)^*)=h(F^*)^*$ for any
$F\in Fin_{or}$ and $\t(E)=E$ for any $E\in Af_{or}$ that is not
of the form $h(F)$ or $h(F)^*$ with $F\in Fin_{or}$.
Note that $\t$ interchanges $E_n^{CC},E_n^{DB}$, interchanges
$E_n^{BB},E_n^{DC}$ and is the identity on all other objects.

We note that $Af_{or}$ is in bijection with the set of (possibly
twisted) affine Cartan matrices in \cite{K83, p.44,45}.

\subhead 1.7\endsubhead
Let $\SS$ be the collection consisting of all (isomorphism classes of)
pairs $(W,\o)$ with $W$ as in 1.4 and $\o\in\Om_W$. Let $(W,\o)\in\SS$.
We say that $\o$ is {\it permissible} if $S^\o\cap S_*\ne\emp$. Let
$\Om^{per}=\{\o\in\Om;\o\text{ permissible}\}$.
Let $\SS^{per}$ be the collection consisting of all (isomorphism classes
of) pairs $(W,\o)\in\SS$ with $\o\in\Om_W^{per}$.

We define a map $\e:\SS^{per}@>>>Af_{or}$, $(W,\o)\m E$ as follows.
If $W$ is simply laced we set $E=W^\o$ (note that $\o$ is good in this
case); if $W$ is not simply laced, so that $\o=1$, $E$ is $W$ with its
canonical orientation. The following result can be deduced from
the tables in \cite{K83, p.44,45}.

(a) $\e:\SS^{per}@>>>Af_{or}$ is a bijection.
\nl
Assuming that $W$ is simply laced we define
$\Om^{adm}=\{\o\in\Om;\o\text{ admissible}\}$ (see 1.2).
Let $\SS^{adm}$ be the collection consisting of all (isomorphism classes
of) pairs $(W,\o)\in\SS$ with $W$ simply laced and $\o\in\Om_W^{adm}$.

We define a map $\e':\SS^{adm}@>>>Af_{or}$ by $(W,\o)\m W^\o$ (note that
$\o$ is good in this case). The following result appears in
\cite{L93, 14.1.4, 14.1.5}.

(b) $\e':\SS^{adm}@>>>Af_{or}$ is a bijection.
\nl
From (a),(b) we deduce that we have a bijection
$$\SS^{per}@>>>\SS^{adm}$$ given by the composition $\e'{}\i\t\e$
where $\t$ is as in 1.6.

We now make a table of all objects of $\SS^{per},\SS^{adm}$.
Each row of this table is of the form

$(W,\o)....E.....E'.....(W',\o')....n\ge?$

where $(W,\o)\in\SS^{per},(W',\o')\in\SS^{adm}$ and
$E=\e(W,\o)\in Af_{or}$, \lb
$E'=\e'(W',\o')\in Af_{or}$, $E'=\t(E)$; $n\ge?$ is a condition on $n\in\NN$.

\mpb

                     TABLE

\mpb

$(\ca_n,1)....\ca_n....\ca_n....(\ca_n,1).... n\ge1$   

$(\ca_n,\o_2)....E_{n/2}^{BC}....E_{n/2}^{BC}....(\cd_{n+2},\o_4).... n\ge2,\text{even}$

$(\ca_n,\o_2)....E_{(n+1)/2}^{BB}....E_{(n+1)/2}^{DC}....
(\cd_{(n+3)/2},\o_2)....n\ge3,\text{odd}$      

$(\cb_n,1)....E_n^{DB}....E_n^{CC}....(\cd_{n+2},\o_{2,2})....
n\ge3$ 

$(\cc_n,1)....E_n^{CC}....E_n^{DB}....(\cd_{2n},\o'_{2,2})....n\ge2$   

$(\cd_n,1)....\cd_n....\cd_n....(\cd_n,1)....n\ge4$
                     
$(\cd_n,\o_2)....E_{n-1}^{DC}....E_{n-1}^{BB}....(\ca_{2n-3},\o_2)....n\ge4$

$(\cd_4,\o_3)....\cg_2^* ....\cg_2^*....(\cd_4,\o_3)$

$(\ce_6,1)....\ce_6 ....\ce_6....(\ce_6,1)$

$(\ce_6,\o_2)....\cf_4^*....\cf_4^*....(\ce_6,\o_2)$ 

$(\ce_7,1)....\ce_7 ....\ce_7....(\ce_7,1)$

$(\ce_8,1)....\ce_8 ....\ce_8....(\ce_8,1)$

$(\cf_4,1)....\cf_4 ....\cf_4....(\ce_7,\o_2)$      

$(\cg_2,1)....\cg_2....\cg_2....(\ce_6,\o_3)$

\head 2. Sharpness\endhead
\subhead 2.1\endsubhead
We now assume that $W,S$ is a Weyl group.
Let $\Irr(W)$ be the set of isomorphism classes of
irreducible representations (over $\bbq$) of $W$.
For any $E\in\Irr(W)$ we denote by $D_{E,q}\in\QQ[q]$ the generic
degree of the
Iwahori-Hecke algebra representation associated to $E$; here $q$
is an indeterminate. Let $N(E)$ be the smallest integer $\ge1$ such that
$N(E)D_{E,q}\in\ZZ[q]$. Let $z(E)$ be the largest integer $n\ge0$ such
that
$D_{E,q}(q+1)^{-n}\in\QQ[q]$. From the known formulas for $D_{E,q}$ one
can see that $z(E)\le r(op)$, (see 1.3).

\subhead 2.2\endsubhead
Let $\o\in\fA_W$. Assuming that $W$ is $\{1\}$ or irreducible, we say
that $(W,\o)$ is {\it sharp} if there exists a special representation
$E\in\Irr(W)$ such that $z(E)=r(op)$ (it is necessarily unique) and if
at least one of (a),(b) below holds.

(a) $\o=op$;

(b) $N(E)$ is divisible by $3$ or $ord(\o)=3$.
\nl
(This definition is similar to, but different from, that in \cite{L20}.)
We now give a list of the various sharp $(W,\o)$. (Here $t\in\NN$.)

(i) $(A_{(t^2-3^2)/8},op)$, $t=1\mod2$;

(ii) $(B_{(t^2-1)/4},op)$, $t=1\mod2$;

(iii) $(D_{t^2/4},op)$, $t=0\mod2$;

(iv) $(D_4,\o_3)$;

(v)  $(E_6,1)$; 

(vi) $(E_8,op),(E_7,op),(E_6,op),(F_4,op),(G_2,op)$.
\nl
In this paper we use the convention:

$D_0=D_1=\{1\}$, $D_2=A_1\T A_1$, $D_3=A_3$, $B_0=\{1\}$, $B_1=A_1$,
$A_0=A_{-1}=\{1\}$.

\subhead 2.3\endsubhead
Let $\o\in\fA_W$. Assuming that $W$ is $\{1\}$ or $\o$-irreducible,
we say that $(W,\o)$ is {\it sharp} if for any irreducible factor
$W'$ of $W$ the following holds:

(a) $(W',\o^n)$ is sharp (see 2.2) where $n$ is the smallest
integer $\ge1$ such that $W'$ is $\o^n$-stable.

For a general $W$ and $\o\in\fA_W$ we say that $(W,\o)$ is {\it sharp}
if $(W',\o)$ is sharp for any $\o$-irreducible factor $W'$ of $W$.

\subhead 2.4\endsubhead
We now assume that $W,S$ is an irreducible affine Weyl group. Let
$\o\in\Om,I\sneq S$. We say that the triple $(W,\o,I)$ is {\it sharp}
if (a),(b) below hold.

(a) $I$ is $\o$-stable and $S-I$ is contained in a single
$\Om$-orbit on $S$.

(b) $(W_I,\o)$ is sharp, see 2.3. 
\nl
Let $\SS_{sharp}$ be the set of (isomorphism classes of) pairs
$(W,\o)$ such that $(W,\o,I)$ is a sharp triple for some $I\sneq S$.

\subhead 2.5\endsubhead
Let $(\QQ_{\ge0}\T\QQ_{\ge0})/\si$ (resp. $(\NN\T\NN))/\si$) be the
set of all unordered pairs of numbers in $\QQ_{\ge0}$ (resp. $\NN$).
We have maps $\io,\io'$ from $(\QQ_{\ge0}\T\QQ_{\ge0})/\si$ to
itself given by
$$\io(a,b)=(a+b,|a-b|), \io'(a,b)=((a+b)/2,|a-b|/2).$$
Note that $\io\io'=\io'\io=1$ hence $\io,\io'$ are bijections.

Let $(W,\o)\in\SS_{sharp}$. Let $\Th=\Th_{W,\o}$ be the set of all
subgroups $W_I$ (up to $\Om^1$-conjugation) of $W$ with $I\sneq S$ such
that $(W,\o,I)$ is sharp.

In the table below
we describe $\Th$ in each case; $W_I$ is specified by its type.
In cases $(b_1)-(g_2)$, $\Th$ is parametrized by a certain subset
$X\sub(\NN\T\NN)/\si$; in each of these cases we also specify the
image $Y$ of $X$ under $\io'$ (when $W=\ca_n$) or under $\io$
(when $W$ is $\cb_n,\cc_n,\cd_n$). We have $Y\sub(\NN\T\NN)/\si$;

\mpb

            TABLE

\mpb

(a)  $(\ca_n,\o_{n+1})$, $\Th=\{1\}$, $n\ge1$.

$(b_1)$ $(\ca_n,\o_2)$, $n\ge2$ even such that
$$X=\{(t,r);t=\pm1\mod8,r=\pm3\mod8,t^2+r^2=8n+10\}\ne\emp,$$
$$\Th=\{W_I=A_{(t^2-9)/8}\T A_{(r^2-9)/8};(t,r)\in X\},$$
$$Y=\{(x,y);x=2\mod4,y=1\mod2,x^2+y^2=4n+5\}.$$                    

$(b_2)$ $(\ca_n,\o_2)$, $n\ge3$ odd such that 
$$X=\{(t,r);t=\pm1\mod8,r=\pm1\mod8,t^2+r^2=8n+10\}\ne\emp,$$
$$\Th=\{W_I=A_{(t^2-9)/8}\T A_{(r^2-9)/8};(t,r)\in X\},$$
$$Y=\{(x,y);x=0\mod4,y=1\mod2,x+y=\pm1\mod8,x^2+y^2=4n+5\}.$$

$(b_3)$ $(\ca_n,\o'_2)$, $n\ge3$ odd such that
$$X=\{(t,r);t=\pm3\mod8,r=\pm3\mod8,t^2+r^2=8n+10\}\ne\emp,$$
$$\Th=\{W_I=A_{(t^2-9)/8}\T A_{(r^2-9)/8};(t,r)\in X\},$$
$$Y=\{(x,y);x=0\mod4,y=1\mod2,x+y=\pm3\mod8,x^2+y^2=4n+5\}.$$

$(c_1)$ $(\cb_n,1)$, $n\ge3$ such that
$$X=\{(t,r);t=0\mod4,r=1\mod2,t^2+r^2=4n+1\}\ne\emp,$$
$$\Th=\{W_I=D_{t^2/4}\T B_{(r^2-1)/4};(t,r)\in X\},$$
$$Y=\{(x,y);x=1\mod2,y=1\mod2,x=\pm y\mod8,x^2+y^2=8n+2\}.$$

$(c_2)$ $(\cb_n,\o_2)$, $n\ge3$ such that
$$X=\{(t,r);t=2\mod4,r=1\mod2,t^2+r^2=4n+1\}\ne\emp,$$
$$\Th=\{W_I=D_{t^2/4}\T B_{(r^2-1)/4};(t,r)\in X\},$$
$$Y=\{(x,y);x=1\mod2,y=1\mod2,x=\pm y+4\mod8,x^2+y^2=8n+2\}.$$

$(d)$ $(\cc_n,1)$, $n\ge2$ such that
$$X=\{(t,r);t=1\mod2,r=1\mod2,t^2+r^2=4n+2\}\ne\emp,$$
$$\Th=\{W_I=B_{(t^2-1)/4}\T B_{(r^2-1)/4};(t,r)\in X\},$$
$$Y=\{(x,y);x=0\mod4,y=2\mod4,x^2+y^2=8n+4\}.$$                

$(e)$ $(\cc_n,\o_2)$, $n\ge2$ such that
$$X=\{(t,r);t=2\mod4,r=1\mod2,t^2+r^2=8n+5\}\ne\emp,$$
$$\Th=\{W_I=B_{(t^2-2^2)/16}\T B_{(t^2-2^2)/16}\T A_{(r^2-3^2)/8};
(t,r)\in X\},$$
$$Y=\{(x,y);x=\pm1\mod8,y=\pm3\mod8,x^2+y^2=16n+10\}.$$
       
$(f_1)$  $(\cd_n,1)$, $n\ge4$ such that
$$X=\{(t,r);t=0\mod4,r=0\mod4,t^2+r^2=4n\}\ne\emp,$$ 
$$\Th=\{W_I=D_{t^2/4}\T D_{r^2/4};(t,r)\in X\},$$
$$Y=\{(x,y);x=0\mod4,y=0\mod4,x=y\mod8,x^2+y^2=8n\}.$$

$(f_2)$ $(\cd_n,\o_{2,2})$, $n\ge4$ such that
$$X=\{(t,r);t=2\mod4,r=2\mod4,t^2+r^2=4n\}\ne\emp,$$ 
$$\Th=\{W_I=\D_{t^2/4}\T D_{r^2/4};(t,r)\in X\},$$
$$Y=\{(x,y);x=0\mod4,y=0\mod4,x=y+4\mod8,x^2+y^2=8n\}.$$

$(f_3)$ $(\cd_n,\o_2)$, $n\ge4$ such that
$$X=\{(t,r);t=0\mod4,r=2\mod4,t^2+r^2=4n\}\ne\emp,$$ 
$$\Th=\{W_I=D_{t^2/4}\T D_{r^2/4};(t,r)\in X\},$$
$$Y=\{(x,y);x=2\mod4,y=2\mod4,x^2+y^2=8n\}.$$

$(g_1)$  $(\cd_n,\o\in{}'\Om^1)$, $n\ge5$ such that
$$X=\{(t,r);t=0\mod4,r=1\mod2,t/4=n\mod2,t^2+r^2=8n+1\}\ne\emp,$$
$$\Th=\{W_I=D_{t^2/16}\T D_{t^2/16}\T A_{(r^2-3^2)/8};(t,r)\in X\},$$
$$Y=\{(x,y);x=\pm1\mod8,y=\pm1\mod8,x^2+y^2=16n+2\}.$$

$(g_2)$ $(\cd_n,\o\in{}'\Om^2)$, $n\ge5$ such that
$$X=\{(t,r);t=0\mod4,r=1\mod2,t/4=n+1\mod2,t^2+r^2=8n+1\}\ne\emp,$$
$$\Th=\{W_I=D_{t^2/16}\T D_{t^2/16}\T A_{(r^2-3^2)/8};(t,r)\in X\},$$
$$Y=\{(x,y);x=\pm3\mod8,y=\pm3\mod8,x^2+y^2=16n+2\}.$$

$(h)$ $(\cd_4,\o_3)$, $\Th=\{D_4\}$.

$(i_1)$ $(\ce_6,1)$, $\Th=\{E_6\}$.

$(i_2)$ $(\ce_6,\o_3)$, $\Th=\{D_4\}$.

$(i_3)$ $(\ce_6,\o_2)$; $\Th=\{E_6\}$.

$(j_1)$ $(\ce_7,1)$, $\Th=\{E_7\}$.

$(j_2)$ $(\ce_7,\o_2)$, $\Th=\{E_6\}$.

$(k)$ $(\ce_8,1)$, $\Th=\{E_8\}$.

$(l)$ $(\cf_4,1)$, $\Th=\{F_4\}$.

$(m)$ $(\cg_2,1)$, $\Th=\{G_2\}$.

\subhead 2.6\endsubhead
Let $(W,\o,I)$ be a sharp triple. We say that $(W,\o,I)$ is
{\it strictly sharp} if (a),(b) below hold.

(a) $S-I$ is an $\Om^1$-orbit in $S$.

(b) In cases 2.5$(b_1)-(g_2)$,the parameters $t,r$ of of $I$ satisfy
$t-r\in\{-1,0,1\}$ (that is, they are as close as possible).

\subhead 2.7\endsubhead
The following is a table of all strictly sharp triples $(W,\o,I)$.
(We again indicate $I$ by the type of $W_I$.)

\mpb

     TABLE

\mpb

$(a)$ $(\ca_n,\o_{n+1},\emp),n\ge1$;

$(c_1)$
$$(\cb_n,1,D_{((x-1)/2)^2/4}\T B_{((x+1)/2)^2-1)/4}),$$
$x\in1+8\NN$, $(x^2-1)/8=n\ge3$, or
$$(\cb_n,1,D_{((x+1)/2)^2/4}\T B_{((x-1)/2)^2-1)/4}),$$
$x\in7+8\NN$, $(x^2-1)/8=n\ge3$;

$(c_2)$
$$(\cb_n,\o_2,D_{((x-1)/2)^2/4}\T B_{((x+1)/2)^2-1)/4}),$$ 
$x\in5+8\NN$,  $(x^2-1)/8=n\ge3$, or
$$(\cb_n,\o_2,D_{((x+1)/2)^2/4}\T B_{((x-1)/2)^2-1)/4}),$$ 
$x\in3+8\NN$,  $(x^2-1)/8=n\ge3$;

$(d)$
$$(\cc_n^a,1,B_{(x^2-4)/16}\T B_{(x^2-4)/16}),$$
$x\in2+4\NN$, $(x^2-4)/8=n\ge2$;

$(e)$
$$(\cc_n,\o_2,
B_{(((x+1)/2)^2-2^2)/16}\T B_{(((x+1)/2)^2-2^2)/16}\T
A_{(((x-1)/2)^2-3^2)/8}),$$
$x\in3+8\NN$, $(x^2-3^2)/16=n\ge2$, or
$$(\cc_n,\o_2, B_{(((x-1)/2)^2-2^2)/16}\T B_{(((x-1)/2)^2-2^2)/16}\T
A_{(((x+1)/2)^2-3^2)/8}),$$
$x\in5+8\NN$, $(x^2-3^2)/16=n\ge2$;

$(f_1)$
$$(\cd_n,1,D_{x^2/16}\T D_{x^2/16}),$$
$x\in8\NN$, $x^2/8=n\ge4$; 

$(f_2)$
$$(\cd_n,\o_{2,2},D_{x^2/16}\T D_{x^2/16}),$$
$x\in4+8\NN$, $x^2/8=n\ge4$;

$(g_1)$
$$(\cd_n,\o\in{}'\Om^1,
D_{(((x+1)/2)^2)/16}\T D_{(((x+1)/2)^2)/16}\T
A_{(((x-1)/2)^2-3^2)/8}),$$
$x\in7+8\NN$, $(x^2-1)/16=n\ge5$, or
$$(\cd_n,\o\in{}'\Om^1,D_{(((x-1)/2)^2)/16}\T D_{(((x-1)/2)^2)/16}\T
A_{(((x+1)/2)^2-3^2)/8}),$$
$x\in1+8\NN$, $(x^2-1)/16=n\ge5$;

$(i_2)$ $(\ce_6,\o_3,D_4)$;

$(j_2)$ $(\ce_7,\o_2,E_6)$;

$(k)$ $(\ce_8,1,E_8)$;

$(l)$ $(\cf_4,1,F_4)$;

$(m)$ $(\cg_2,1,G_2)$.
\nl
We see that in each case we must have $\o\in\Om^1$ and the set $I$
is uniquely determined by $(W,\o)$, if it exists. 

\subhead 2.8\endsubhead
Let $W$ be an irreducible affine Weyl group. Let $W'$ be an
irreducible affine Weyl group such that $W,W'$ with their canonical
orientations are of the form $h(F),h(F^*)$ where $F\in Fin_{or}$.
Let $\o'\in\Om_{W'}^1$.

We say that $(W,\o')$ is a {\it blunt pair} if $(W',\o',I)$ is strictly
sharp for some $I\sneq S'$ (which is necessarily unique, see 2.7). The
following is a table of the blunt pairs $(W,\o')$; it can be obtained
from that of the strictly sharp triples (see 2.7).

\mpb

            TABLE

\mpb

$(a)$ $(\ca_n,\o_{n+1})$, $n\ge1$;

$(c_1)$ $(\cc_n,1)$, $n\ge3$, $x\in(1+8\NN)\cup(7+8\NN)$, $n=(x^2-1)/8$;

$(c_2)$ $(\cc_n,\o_2)$, $n\ge3$, $x\in(3+8\NN)\cup(5+8\NN)$,
$n=(x^2-1)/8$;

$(d)$ $(\cb_n,1)$, $n\ge4$, $x\in2+4\NN$, $n=(x^2-4)/8$;

$(e)$ $(\cb_n,\o_2)$, $n\ge3$, $x\in(3+8\NN)\cup(5+8\NN)$,
$n=(x^2-3^2)/16$;

$(f_1)$ $(\cd_n,1)$, $n\ge4$, $x\in8\NN$, $n=x^2/8$; 

$(f_2)$ $(\cd_n,\o_{2,2})$, $n\ge4$, $x\in4+8\NN$, $n=x^2/8$;

$(g_1)$ $(\cd_n,\o\in{}'\Om_{W'}^1)$, $n\ge5$,
$x\in(1+8\NN)\cup(7+8\NN)$, $n=(x^2-1)/16$;

$(i_2)$ $(\ce_6,\o_3)$;

$(j_2)$ $(\ce_7,\o_2)$;

$(k)$ $(\ce_8,1)$;

$(l)$ $(\cf_4,1)$;

$(m)$ $(\cg_2,1)$.

\head 3. Bluntness\endhead
\subhead 3.1\endsubhead
Assume that $W_1$ is an irreducible oriented Weyl group. We say that
$W_1$ is {\it blunt} if there exists a blunt pair $(W,\o')$ (see 2.8)
such that $h(W_1)=W$ with the canonical orientation ($h$ is
as in 1.6).

The following is a table of the irreducible blunt oriented Weyl groups.

\mpb

             TABLE

\mpb

(a) $A_n$, $n\ge1$;

(b) $B_n$, $n\ge4$, $x\in2+4\NN$, $(x^2-4)/8=n$;

(c) $C_n$, $n\ge3$, $x\in1+2\NN$, $(x^2-1)/8=n$;

(d) $D_n$, $n\ge4$, $x\in4\NN$, $x^2/8=n$;

(b') $B_n$, $n\ge3$, $x\in(5+8\NN)\cup(3+8\NN)$, $(x^2-3^2)/16=n$;

(d') $D_n$, $n\ge5$, $x\in(1+8\NN)\cup(7+8\NN)$, $(x^2-1)/16=n$;   

$E_6,E_7,E_8,F_4,G_2$.
\nl
The blunt oriented Weyl groups in (b),(c),(d),(b'),(d') can be
arranged in a sequence

$D_2,C_3,B_4,D_5,C_6,B_7,D_8,C_{10},B_{12},D_{14},\do$.
\nl
A (not necessarily irreducible) oriented Weyl group $W$ is said to
be blunt if any of its irreducible factors (with the induced
orientation) is blunt.

\subhead 3.2\endsubhead
We now assume that $E'$ is an irreducible affine Weyl group
with set of simple reflections $S'$ and that $|S'|\ge3$.

For any $J\sneq S'$ let $E'_J$ the Weyl subgroup of $E'$
generated by $J$.

We shall define a function $c:S'@>>>\ZZ_{>0}$.
Assume first that $E'$ is simply laced. Following Coxeter we define
$c:S'@>>>\ZZ_{>0}$ as the unique function such that for any $s\in S'$
we have $c(s)=(1/2)\sum_{s'\in S';ss'\ne s's}c(s')$ (harmonicity) and
$c(s'')=1$ for some $s''\in S'$. (We then have $s''\in S'_*$.)

\subhead 3.3\endsubhead
We now assume that $E'$ is not simply laced. There is a well defined
$(\tW,\ti\o)\in\SS^{adm}$ where $\tW$ has simple reflections $\tS$ and
$\ti\o\in\Om^1_{\tW}$ such that $E'=\tW^{\ti\o}$, $S'=\tS_{\ti\o}$
(see 1.1) with the canonical orientation of $E'$.
(If $E'$ is $\cb_n,\cc_n,\cf_4,\cg_2$ then $(\tW,\ti\o)$ is
$$(\cd_{2n},\o'_{2,2}),(\cd_{n+2},\o_{2,2}),(\ce_7,\o_2),(\ce_6,\o_3)$$
respectively, see 1.7.)
We define $c:S'@>>>\ZZ_{>0}$ as the unique function such that the
composition $\tS@>>>\tS_{\ti\o}=S'@>c>>\ZZ_{>0}$ is as in 3.2 (for
$\tW$ instead of $E'$).

\subhead 3.4\endsubhead
Now let $s\in S'$. Then $E'_{S'-\{s\}}$ is a product
$$(W_1)_{S_1-\{s_1\}}\T(W_2)_{S_2-\{s_2\}}\T\do\T
(W_k)_{S_k-\{s_k\}}$$
where $(W_1,S_1),(W_2,S_2),\do,(W_k,S_k)$ are irreducible affine Weyl
groups and $s_1\in(S_1)_*,s_2\in(S_2)_*,\do,s_k\in(S_k)_*$.

For $i=1,\do,k$, let $n_i$ be the maximum of the
orders of various elements in $\Om^1_{W_i}$. Let
$n(s)=\max_{i\in[1,k]}n_i$. We obtain a function $n:S'@>>>\ZZ_{>0}$.

\subhead 3.5\endsubhead
One can verify that $n(s)/c(s)\in\ZZ_{>0}$ for all $W$ and all
$s\in S'$. We now list the fractions $\{n(s)/c(s);s\in S'\}$ 
for various $E'$.

$E'=\ce_6$: 
$$\left(\matrix 3/1 & 6/2 & 3/3 &6/2&3/1\\
                {} & {} & 6/2 & {}& {}\\
                {} & {} & 3/1 & {}& {}\endmatrix\right)$$
$E'=\ce_7$: 
$$\left(\matrix 2/1 & 2/2 & 6/3 &4/4& 6/3&2/2&2/1\\
                {} & {} & {} &  8/2 {}& {}&{}&{}\endmatrix\right)$$

$E'=\ce_8$:
$$\left(\matrix 2/2 & 8/4 & 6/6 &5/5& 4/4&3/3&2/2&1/1\\
                {} & {} &9/3 &{}& {}&{}&{}&{}\endmatrix\right)$$

$E'=\cf_4$: $\{4/4,3/3,2/2,2/2,1/1\}$

$E'=\cg_2$: $\{3/3,2/2,1/1\}$.

For $a\ge1$ let $S'_a$ be a set of representatives for the
$\Om_{E'}^1$-orbits on

$\{s\in S';n(s)/c(s)=a\}$
\nl
and let

$U_a(E')=\{E'_{S'-\{s\}};s\in S'_a\}$.
\nl
Let

$S'_a{}^b=\{s\in S'_a;E'_{S'-\{s\}}\text{ is blunt}\}$
\nl
and let $U_a^b(E')=\{E'_{S'-\{s\}};s\in S'_a{}^b\}$.

If $E'=\ce_6$ we have $U_1(E')=U_1^b(E')=\{A_2\T A_2\T A_2\}$
$U_3(E')=U_3^b(E')=\{E_6,A_5\T A_1\}$.
                               
If $E'=\ce_7$ we have
$U_1(E')=\{D_6\T A_1,A_3\T A_3\T A_1\}$,
$U_1^b(E')=\{A_3\T A_3\T A_1\}$,
$U_2(E')=U_2^b(E')=\{E_7,A_5\T A_2\}$.

If $E'=\ce_8$ we have
$$U_1(E')=U_1^b(E')=\{E_8,E_7\T A_1,E_6\T A_2,D_5\T A_3,A_4\T A_4,
A_5\T A_2\T A_1,D_8\}.$$
            
We now describe $S'_a{}^b$ and the values of $s\m c(s)$ for
$s\in S'_a{}^b$
for $E'$ of type $\ce_6,\ce_7,\ce_8$ and certain $a$.

$E'=\ce_6,a=1$:
$$\left(\matrix - & - & 3 &-&-\\
                {} & {} & - & {}& {}\\
                {} & {} & - & {}& {}\endmatrix\right)$$
$E'=\ce_6,a=3$:
$$\left(\matrix 1 & 2 & - &-&-\\
                {} & {} & - & {}& {}\\
                {} & {} & - & {}& {}\endmatrix\right)$$
$E'=\ce_7,a=1$: 
$$\left(\matrix - & - & - &4& -&-&-\\
                {} & {} & {} &  - {}& {}&{}&{}\endmatrix\right)$$
$E'=\ce_7,a=2$: 
$$\left(\matrix 1 & - & 3 &-& 3&-&1\\
                {} & {} & {} &  - {}& {}&{}&{}\endmatrix\right)$$
$E'=\ce_8,a=1$:
$$\left(\matrix 2 & - & 6 &5& 4&3&2&1\\
                {} & {} &- &{}& {}&{}&{}&{}\endmatrix\right)$$
The last picture can be refined to the picture
$$\left(\matrix g'_2&-&g_6&g_5&g_4&g_3&g_2&1\\
                {}&{}&-&{}&{}&{}&{}&{}\endmatrix\right)$$
which gives an imbedding of the set of conjugacy classes in the
symmetric group $\fS_5$ into $S$. (We use the notation of
\cite{L84a,p.79,80} for the conjugacy classes in $\fS_5$; the order
of an element in a conjugacy class is equal to $c(s)$ for the
corresponding $s\in S$.

There is a similar imbedding of the set of conjugacy classes in the
symmetric group $\fS_4$ or $\fS_5$ into $S$ for $W$ of type $\cf_4$
or $\cg_2$.
$$\left(\matrix g'_2&g_4&g_3&g_2&1\endmatrix\right)$$
$$\left(\matrix g'_2&g_3&1\endmatrix\right)$$

If $E'=\cf_4$ we have
$$U_1(E')=U_1^b(E')=\{F_4,C_3\T A_1,A_2\T A_2,A_3\T A_1,B_4\}.$$
 The corresponding values
of the function $s\m c(s)$ on $S'_1{}^b$ are $1,2,3,4,2$.

If $E'=\cg_2$ we have
$$U_1(E')=U_1^b(E')=\{G_2,A_1\T A_1,A_2\}.$$
 The corresponding values
of the function $s\m c(s)$ on $S'_1{}^b$ are $1,2,3$.

\subhead 3.6\endsubhead
Assume now that $E'$ is not simply laced.
There is a well defined $(W^!,\o^!)\in\SS^{adm}$ where $W^!$ has simple 
reflections $S^!$ and $\o^!\in\Om_{W^!}-\Om^1_{W^!}$ such that
$E'=(W^!)^{\o^!}$, $S=S^!_{\o^!}$ (see 1.1) with the opposite of the
canonical orientation.
(If $E'$ is $\cb_n,\cc_n,\cf_4,\cg_2$ then $(W^!,\o^!)$ is
$$(\cd_{n+1},\o_2),(\ca_{2n-1},\o_2),(\ce_6,\o_2),(\cd_4,\o_3)$$
respectively, see 1.7.)
We define $c^!:S'@>>>\ZZ_{>0}$ as the unique function such that the
composition $S^!@>>>S^!_{\o^!}=S'@>c^!>>\ZZ_{>0}$ is as in 3.2
(for $W^!$ instead of $E'$).
One can verify that $n(s)/c^!(s)\in\ZZ_{>0}$ for all $s\in S'$.

We now list the fractions $\{n(s)/c^!(s);s\in S'\}$ for $E'$ of certain
types.

$E'=\cf_4$: $\{1/1,2/2,3/3,2/1,4/2\}$

$E'=\cg_2$: $\{1/1,2/2,3/1\}$.

For $a\ge1$ let $S'{}^!_a$ be a set of representatives for the
$\Om_{E'}^1$-orbits on $\{s\in S';n(s)/c^!(s)=a\}$ and let
$$U_a^!(E')=\{E'_{S'-\{s\}};s\in S'_a{}^!\}.$$
Let $S'_a{}^{!b}=\{s\in S'_a{}^!;E'_{S'-\{s\}}\text{ is blunt}\}$
and let
$U_a^{!b}(E')=\{E'_{S'-\{s\}};s\in S'_a{}^{!b}\}$.

If $E'=\cf_4$ we have $U^!_1(E')=\{F_4,B_3\T A_1,A_2\T A_2\}$,
$U^{!b}_1(E')=\{F_4,A_2\T A_2\}$. The corresponding values
of the function $s\m c^!(s)$ on $S'_1{}^{!b}$ are $1,3$.
                    
If $E'=\cg_2$ we have $U^!_1(E')=U^{!b}_1(E')=\{G_2,A_1\T A_1\}$.
 The corresponding values
of the function $s\m c^!(s)$ on $S'_1{}^{!b}$ are $1,2$.

\subhead 3.7\endsubhead
For any $m\ge1$ let $\mu_m$ be the set of primitive $m$-th roots of
$1$ in $\bbq^*$.
For $E'$ of type $\ce_6,\ce_7,\ce_8,\cf_4,\cg_2$ and $a\ge1$ let
$\tS_a^b(W)$ be the multiset whose objects are the elements
$s\in S'_a{}^b(E')$ with $s$ repeated $|\mu_{c(s)}|$ times.
Let $\ti U_a^b(E')$ be the multiset
$\{E'_{S'-\{s\}};s\in \tS_a^b(E')\}$; this is in obvious
bijection with $\tS_a^b(E')$.

For $W$ of type $\cf_4,\cg_2$ let
$\tS_1^{!b}(E')$ be the multiset whose objects are the elements
$s\in S'_1{}^{!b}(E')$ with $s$ repeated $|\mu_{c^!(s)}|$ times.
Let $\ti U_1^{!b}(E')$ be the multiset
$\{E'_{S'-\{s\}};s\in\tS_1^{!b}(E')\}$; this is in obvious
bijection with $\tS_1^{!b}(E')$.

\subhead 3.8\endsubhead
We now assume that $E'$ in 3.2 is associated to an element
$(W,\o)\in\SS_{sharp}$ as follows. There is a well defined
$\o_0\in\Om^{per}$
 (up to conjugation by $\Om^1$) such that $\o_0\in\o\Om^1$.
Consider the row of the table in 1.7
with first item $(W,\o_0)$. The third item in that row
is an object $E'=E'_{W,\o}\in Af_{or}$.

If $(W,\o)=(\cd_4,\o_3)$, so that $E'=\cg_2^*$, we set
$\tS_{W,\o}=\tS_1^{!b}(\cg_2)$,
$\tU_{W,\o}=\tU_1^{!b}(\cg_2)$.

If $(W,\o)$ is one of $(\ce_6,1),(\ce_7,1),(\ce_8,1),(\cf_4,1),
(\cg_2,1)$ so that $E'=W$, we set
$\tS_{W,\o}=\tS_1^b(W)$, $\tU_{W,\o}=\tU_1^b(W)$.

If $(W,\o)=(\ce_6,\o_3)$, so that $E'=\ce_6$, we set
$\tS_{W,\o}=\tS_3^b(\ce_6)$, $\tU_{W,\o}=\tU_3^b(\ce_6)$.

If $(W,\o)=(\ce_6,\o_2)$ so that $E'=\cf_4^*$ we set
$\tS_{W,\o}=\tS_1^{!b}(\cf_4)$, $\tU_{W,\o}=\tU_1^{!b}(\cf_4)$.

If $(W,\o)=(\ce_7,\o_2)$, so that $E'=\ce_7$, we set
$\tS_{W,\o}=\tS_2^b(\ce_7)$, $\tU_{W,\o}=\tU_2^b(\ce_7)$.

Note that in each of these cases there is an obvious bijection
$\b:\tS_{W,\o}@>\si>>\tU_{W,\o}$.

\subhead 3.9\endsubhead
Let $(W,\o)\in\SS_{sharp}$. We associate $E'$ to $(W,\o)$ as in 3.8.

Below we define a nonempty set $\Si=\Si_{W,\o}$ of blunt Weyl
subgroups $E'_J$ of $E'$ (with induced orientation) where
$J\sneq S',|J|=|S'|-1$ (modulo the conjugation action of $\Om_{E'}$).

$(a)$ If $(W,\o)=(\ca_n,\o_{n+1}),n\ge1$, then $\Si=\{A_n\}$.

$(b_1)$ If $(W,\o)=(\ca_n,\o_2)$ (as in 2.5$(b_1)$), so that
$E'=W_{n/2}^{BC}$, then
$$\Si=\{E'_J=B_{(x^2-4)/8}\T C_{(y^2-1)/8};(x,y)\in Y\},$$
with $Y$ as in 2.5$(b_1)$.

$(b_2)$ If $(W,\o)=(\ca_n,\o_2)$ (as in 2.5$(b_2)$), so that
$E'=W_{(n+1)/2}^{DC}$, then
$$\Si=\{E'_J=D_{x^2/8}\T C_{(y^2-1)/8};(x,y)\in Y\},$$
with $Y$ as in 2.5$(b_2)$.  

$(b_3)$ If $(W,\o)=(\ca_n,\o'_2)$ (as in 2.5$(b_3)$), so that
$E'=W_{(n+1)/2}^{DC}$, then
$$\Si=\{E'_J=D_{x^2/8}\T C_{(y^2-1)/8};(x,y)\in Y\},$$
with $Y$ as in 2.5$(b_3)$.

$(c_1)$ If $(W,\o)=(\cb_n,1)$ (as in 2.5$(c_1)$), so that $E'=E_n^{CC}$,
then
$$\Si=\{E'_J=C_{(x^2-1)/8}\T C_{(y^2-1)/8};(x,y)\in Y\},$$
with $Y$ as in 2.5$(c_1)$.

$(c_2)$ If $(W,\o)=(\cb_n,\o_2)$ (as in 2.5$(c_2)$), so that
$E'=E_n^{CC}$, then
$$\Si=\{E'_J=C_{(x^2-1)/8}\T C_{(y^2-1)/8};(x,y)\in Y\},$$
with $Y$ as in 2.5$(c_2)$.

$(d)$ If $(W,\o)=(\cc_n,1)$ (as in 2.5$(d)$), so that $E'=E_n^{DB}$, then
$$\Si=\{E'_J=D_{x^2/8}\T B_{(y^2-4)/8};(x,y)\in Y\},$$
with $Y$ as in 2.5$(d)$.

$(e)$ If $(W,\o)=(\cc_n,\o_2)$ (as in 2.5$(e)$), so that $E'=E_n^{DB}$,
then
$$\Si=\{E'_J=D_{(x^2-1)/16}\T B_{(y^2-9)/16};(x,y)\in Y\},$$
with $Y$ as in 2.5$(e)$.

$(f_1)$ If  $(W,\o)=(\cd_n,1)$ (as in 2.5$(f_1)$), so that $E'=\cd_n$,
then
$$\Si=\{E'_J=D_{x^2/8}\T D_{y^2/8};(x,y)\in Y\},$$
with $Y$ as in 2.5$(f_1)$.

$(f_2)$ If  $(W,\o)=(\cd_n,\o_{2,2})$ (as in 2.5$(f_2)$), so that
$E'=\cd_n$, then
$$\Si=\{E'_J=D_{x^2/8}\T D_{y^2/8};(x,y)\in Y\},$$
with $Y$ as in 2.5$(f_2)$.

$(f_3)$ If  $(W,\o)=(\cd_n,\o_2)$ (as in 2.5$(f_3)$), so that
$E'=W_{n-1}^{BB}$, then
$$\Si=\{E'_J=B_{(x^2-4)/8}\T B_{(y^2-4)/8};(x,y)\in Y\},$$
with $Y$ as in 2.5$(f_3)$.

$(g_1)$ If $(W,\o)=(\cd_n,\o\in{}'\Om^1)$ (as in 2.5$(g_1)$), so that
$E'=\cd_n$, then
$$\Si=\{E'_J=D_{(x^2-1)/16}\T D_{(y^2-1)/16};(x,y)\in Y\},$$
with $Y$ as in 2.5$(g_1)$.

$(g_2)$ If $(W,\o)=(\cd_n,\o\in{}'\Om^2)$ (as in 2.5$(g_2)$), so that
$E'=W_{n-1}^{BB}$, then
$$\Si=\{E'_J=B_{(x^2-9)/16}\T B_{(y^2-9)/16};(x,y)\in Y\},$$
with $Y$ as in 2.5$(g_2)$.

$(h)$ If $(W,\o)=(\cd_4,\o_3)$, so that $E'=\cg_2^*$, then
$$\Si=\tU_{W,\o}=\{G_2,A_1\T A_1\},$$
(see 3.8).     

$(i_1)$ If $(W,\o)=(\ce_6,1)$, so that $E'=\ce_6$, then
$$\Si=\tU_{W,\o}=\{A_2\T A_2\T A_2,A_2\T A_2\T A_2\},$$
(see 3.8).

$(i_2)$ If $(W,\o)=(\ce_6,\o_3)$, so that $E'=\ce_6$, then
$$\Si=\tU_{W,\o}=\{E_6,A_5\T A_1\},$$
(see 3.8).

$(i_3)$ If $(W,\o)=(\ce_6,\o_2)$, so that $E'=\cf_4^*$, then
$$\Si=\tU_{W,\o}=\{F_4,A_2\T A_2,A_2\T A_2\},$$
(see 3.8).

$(j_1)$ If $(W,\o)=(\ce_7,1)$, so that $E'=\ce_7$, then
$$\Si_=\tU_{W,\o}=\{A_3\T A_3\T A_1,A_3\T A_3\T A_1\},$$
(see 3.8).

$(j_2)$ If $(W,\o)=(\ce_7,\o_2)$, so that $E'=\ce_7$, then
$$\Si=\tU_{W,\o}=\{E_7,A_5\T A_2,A_5\T A_2\},$$
(see 3.8).

$(k)$ If $(W,\o)=(\ce_8,1)$, so that $E'=\ce_8$, then
$$\align&\Si=\tU_{W,\o}=
\{E_8,E_7\T A_1,E_6\T A_2,E_6\T A_2, D_5\T A_3,D_5\T A_3,\\&
A_4\T A_4,A_4\T A_4,A_4\T A_4,A_4\T A_4,A_5\T A_2\T A_1,
A_5\T A_2\T A_1,D_8\},\endalign$$
(see 3.8).

$(l)$ If $(W,\o)=(\cf_4,1)$, so that $E'=\cf_4$,  then
$$\Si=\tU_{W,\o}=
\{F_4,C_3\T A_1,A_2\T A_2,A_2\T A_2,A_3\T A_1,A_3\T A_1,B_4\},$$
(see 3.8).

$(m)$ If $(W,\o)=(\cg_2,1)$, so that $E'=\cg_2$, then
$$\Si=\tU_{W,\o}=\{G_2,A_1\T A_1,A_2,A_2\},$$
(see 3.8).
\nl
In cases $(b_1)-(g_2)$ let $\a:\Th_{W,\o}@>si>>\Si_{W,\o}$ be the
bijection induced by the bijection $X@>\si>>Y$ (see 2.5) defined by
$\io'$ in 2.5 (when $W=\ca_n$) or by $\io$ in 2.5
(when $W$ is $\cb_n,\cc_n,cd_n$). 

\subhead 3.10\endsubhead
Let $\ti\SS_{sharp}$ be the multiset of all $((W,\o),W_I,*)$
where $(W,\o)\in\SS_{sharp}$ and:

(in cases 2.5$(a)-(g_2)$), $W_I\in\Th_{W,\o}$ and $*=\{\}$;

(in cases 2.5$(h)-(m)$), $W_I$ is the unique object of $\Th_{W,\o}$
and $*$ runs through $\tS_{W,\o}$ (see 3.8).

Let $\ti\SS_{blunt}$ be the multiset of all $((W,\o),E'_J)$ where
$(W,\o)\in\SS_{sharp}$ and $E'_J\in\Si_{W,\o}$.

There is a unique bijection $\ti\SS_{sharp}@>\si>>\ti\SS_{blunt}$
given by

$((W,\o),W_I,*)\m((W,\o),\a(W_I))$
\nl
if $(W,\o)$ is as in 2.5$(a)-(g_2)$ and

$((W,\o),W_I,*)\m((W,\o),\b(*))$
\nl
if $(W,\o)$ is as
in 2.5$(h)-(m)$; here $\a:\Th_{W,\o}@>>>\Si_{W,\o}$ is the obvious
bijection in case 2.5$(a)$ and is the bijection in 3.9 in cases
2.5$(b_1)-(g_2)$; $\b$ is the bijection in 3.8.

\head 4. Comments\endhead
\subhead 4.1\endsubhead
The sharp pairs $(W,\o)$ listed in 2.2 are exactly the cases where
$Cu(W,\o)\ne\emp$ (see 0.1).

\subhead 4.2\endsubhead
The classification of blunt pairs given in 2.8
is the same as the classification (given in
\cite{L84b, p.206}) of the cuspidal local systems on unipotent classes
in the simply connected almost simple group $G$ over $\CC$ corresponding
to the affine Weyl group $W'$. Then $\o'$ can be viewed as the character
by which the centre of $G$ acts on the cuspidal local system.

\subhead 4.3\endsubhead
We now assume that $W$ is an irreducible Weyl group and $\o\in\fA_W$.
We can find an irreducible affine Weyl group
$\tW,\tS$, $\ti\o\in\Om_{\tS}$ and $s\in\tS_*$ such that
$\ti\o(s)=s$, $W=\tW_{\tS-\{s\}}$ and $\o=\ti\o|_W$.
Assume now that $(W,\o)$ is one of

$(D_4,\o_3),(E_6,1),(E_6,\o_2),(E_7,1),(E_8,1),(F_4,1),(G_2,1)$
\nl
so that $(\tW,\ti\o)$ is 

$(\cd_4,\o_3),(\ce_6,1),(\ce_6,\o_2), (\ce_7,1),(\ce_8,1),
(\cf_4,1),(\cd_2,1)$
\nl
respectively.
Then $|\Si_{\tW,\ti\o}|$ is $2,2,3,2,13,7,4$ respectively.
We see that $|\Si_{\tW,\ti\o}|=|Cu(W,\o)|$, see 0.1.

\subhead 4.4\endsubhead
The bijection $\tS_{sharp}@>\si>>\tS_{blunt}$ in 3.10 is essentially
the same as the bijection between the arithmetic diagrams and
geometric diagrams in \cite{L95},\cite{L02} (except that certain
repetitions are ignored).

\subhead 4.5\endsubhead
In this subsection we define an involution
$\z_{fin}:\SS^{adm}_{fin}@>>>\SS^{adm}_{fin}$
(resp. an involution $\z:\SS^{adm}@>>>\SS^{adm}$).

Let $(W,\o)$ be in $\SS^{adm}_{fin}$ (resp. in $\SS^{adm}$).
Let $S$ be the set of simple reflections in $W$. We want to define a
new object $(\hW,\hat\o)$ in $\SS^{adm}_{fin}$ (resp. in $\SS^{adm}$).
If $\o=1$ we set $(\hW,\hat\o)=(W,1)$. We now assume that $\o\ne1$.
Let $m\ge2$ be the order of $\o$. Let $\ce$ be the set of $2$-element
subsets $\{s,s'\}$ of $S$ such that $ss'\ne s's$. For any $x\in S_\o$
let $\hS_x$ be a set with a transitive
$\ZZ/m$-action such that $|\hS_x|=m/|x|$. Let
$\hS=\sqc_{x\in S_\o}\hS_x$ with the $\ZZ/m$-action which for any
$x\in S_\o$ restricts to the $\ZZ/m$-action on $\hS_x$ considered
above. Let $\hat\o:\hS@>>>\hS$ be the action of the generator
$1\in\ZZ/m$.
For any unordered pair $x,x'$ of distinct elements of $S_\o$
let $\ce_{x,x'}$ be the set of elements of $\ce$ of the
form $\{s,s'\}$ for some $s\in x,s'\in x'$. We define a set
$\hat\ce_{x,x'}$ of two-element subsets of
$\hS$ as follows. If $\ce_{x,x'}=\emp$ then $\hat\ce_{x,x'}=\emp$.
If $\ce_{x,x'}\ne\emp$ we choose $\x\in\hS_x,\x'\in\hS_{x'}$ and we
define $\hat\ce_{x,x'}$ to be the orbit of $\{\x,\x'\}$
for the $\ZZ/m$-action on the set of two-element subsets
of $\hS$ induced by the $\ZZ/m$-action on $\hS$.
Let $\hat\ce$ be the (disjoint) union of the sets $\hat\ce_{x,x'}$ for
various $x,x'$ as above.
Now $\hW$ is the Coxeter group for which the Coxeter graph has vertices
$\hS$ and edges $\hat\ce$. Now $\hat\o$ can be regarded as an
automorphism of $\hW$ preserving $\hS$. One can verify that
$(\hW,\hat\o)$ is in $\SS^{adm}_{fin}$ (resp. in $\SS^{adm}$) and
that it is independent of
the choices. From the definition we see that the map
$\z_{fin}:\SS^{adm}_{fin}@>>>\SS^{adm}_{fin}$
(resp. $\z:\SS^{adm}@>>>\SS^{adm}$) given by $(W,\o)\m(\hW,\hat\o)$ 
corresponds under the bijection
$\SS^{adm}_{fin}@>>>Fin_{or}$ in 1.3(a)
(resp. $\SS^{adm}@>>>Af_{or}$ in 1.7(b)) 
to the involution $F\m F^*$ (see 1.1) of $Fin_{or}$ (resp. $Af_{or}$).
In particular $\z_{fin},\z$ are involutions.

The involution $\z_{fin}:\SS^{adm}_{fin}@>>>\SS^{adm}_{fin}$
interchanges

$(A_{2n-1},\o_2)$ with $(D_{n+1},\o_2)$ for $n\ge3$;
\nl
it is the identity on all other objects.

The involution $\z:\SS^{adm}@>>>\SS^{adm}$ interchanges

$(\cd_{n+1},\o_2)$ with $(\cd_{2n},\o'_{2,2})$ for $n\ge2$;

$(\cd_{n+2},\o_{2,2})$ with $(\ca_{2n-1},\o_2)$ for $n\ge3$;

$(\cd_4,\o_3)$ with $(\ce_6,\o_3)$;

$(\ce_6,\o_2)$ with $(\ce_7,\o_2)$;
\nl
it is the identity on all other objects.

\subhead 4.6\endsubhead
As mentioned in \cite{L02, 6.26}, the discussion in
\cite{L02,\S6} has been influenced by the paper \cite{K69} on outer
automorphisms of finite order of a complex simple Lie algebra. (At
the time of writing \cite{L02} I was not aware of the paper \cite{S56}
where a similar study was made earlier than that in \cite{K69}.)


\widestnumber\key{KL87}
\Refs
\ref\key{S56}\by J.de Siebenthal\paper Sur les groupes compactes non
connexes\jour Comment. Math. Helv.\vol31\yr1956\pages41-89\endref
\ref\key{K69}\by V.G.Kac\paper Automorphisms of finite order of
semisimple Lie algebras\jour Funkt.Analiz i ego prilozh.\vol3\yr1969
\pages94-96\endref
\ref\key{K83}\by V.G.Kac\book Infinite dimensional Lie algebras\publ
Birkh\"auser\yr1983\endref
\ref\key{KL87}\by D.Kazhdan and G.Lusztig\paper Proof of the
Deligne-Langlands conjecture for Hecke \lb
algebras\jour Inv. Math.\vol87
\yr1987\pages153-215\endref
\ref\key{LS79}\by G.Lusztig and N.Spaltenstein\paper Induced unipotent
classes\jour J. Lond. Math. Soc.\vol19\yr1979\pages41-52\endref
\ref\key{L83}\by G.Lusztig\paper Some examples of square integrable
representations of semisimple p-adic \lb groups
\jour Trans. Amer. Math. Soc.
\vol227\yr1983\pages623-653\endref
\ref\key{L84a}\by G.Lusztig\book Characters of reductive groups over a
finite field\bookinfo Ann. Math. Studies 107\publ Princeton U. Press\yr
1984\endref
\ref\key{L84b}\by G.Lusztig\paper Intersection cohomology complexes
on a reductive group\jour Invent. Math.\vol75\yr1984\pages205-272\endref 
\ref\key{L93}\by G.Lusztig\book Introduction to quantum groups\bookinfo
Progr. in Math. 110\publ Birkh\"auser Boston\yr1993\endref
\ref\key{L95}\by G.Lusztig\paper Classification of unipotent
representations of simple $p$-adic groups\jour Int. Math. Res. Notices
\yr1995\pages517-589\endref
\ref\key{L02}\by G.Lusztig\paper Classification of unipotent
representations of simple $p$-adic groups, II\jour Represent. Th.\vol6
\yr2002\pages243-289\endref
\ref\key{L14}\by G.Lusztig\book Hecke algebras with unequal parameters
\bookinfo version 2,  arxiv:math/0208154\endref
\ref\key{L20}\by G.Lusztig\paper Unipotent blocks and weighted affine
Weyl groups\jour arxiv:2010.02095\endref
\endRefs
\enddocument